\documentclass[11pt]{article}

\usepackage{amssymb, amsmath, amsthm}

\usepackage{epsfig}
\usepackage{graphicx}
\usepackage{color}
\usepackage{mathptmx}

\usepackage{hyperref}
\usepackage{fullpage}
\usepackage{pinlabel}
\theoremstyle{plain}
\newtheorem{theorem}{Theorem}

\newtheorem{question}[theorem]{Question}
\newtheorem{conjecture}[theorem]{Conjecture}

\newcommand{\R}{{\mathbb R}}

\newcommand{\defn}[1]{\textbf{#1}}
\newcommand{\boundary}{\partial}

\newcommand{\vol}{\operatorname{vol}}

\newcommand{\spacing}{ \parskip 6.6pt \parindent 0pt}

\spacing

\begin{document}

\title{Abstractly Planar Spatial Graphs} 
\author{Scott A. Taylor}

\maketitle

\begin{abstract}
This is a survey article for the forthcoming ``A Concise Encyclopedia of Knot Theory.'' We focus on the topology of spatial graphs with few vertices and edges, paying particular attention to Brunnian $\theta$-graphs.
\end{abstract}

Beginning courses in graph theory prove many wonderful theorems about planar graphs. An even more wonderful theory arises when we put planar graphs (which we'll henceforth refer to as \defn{abstractly planar} graphs) into 3-dimensional space. One way of doing this is to choose an embedding of an abstractly planar graph $G$ in the sphere $S^2$ and then include $S^2$ into the 3-sphere $S^3 = \R^3 \cup \{\infty\}$ (tamely). As an abstractly planar graph may have several planar embeddings, we may wonder if we can end up with different embeddings in $S^3$. It turns out that in $S^3$ all planar embeddings give rise to \defn{equivalent} (that is, ambient isotopic) spatial graphs \cite{Mason}. Such spatial graphs are called \defn{trivial}.  A nontrivial spatial graph is  \defn{knotted}.  Are there knotted embeddings? Yes there are! 

\begin{figure}
\labellist
\small\hair 2pt
\pinlabel $\centerdot \centerdot \centerdot$ at 596 129
\pinlabel \rotatebox{45}{$\centerdot \centerdot \centerdot$} at 711 142
\endlabellist
\centering
\includegraphics[scale=0.35]{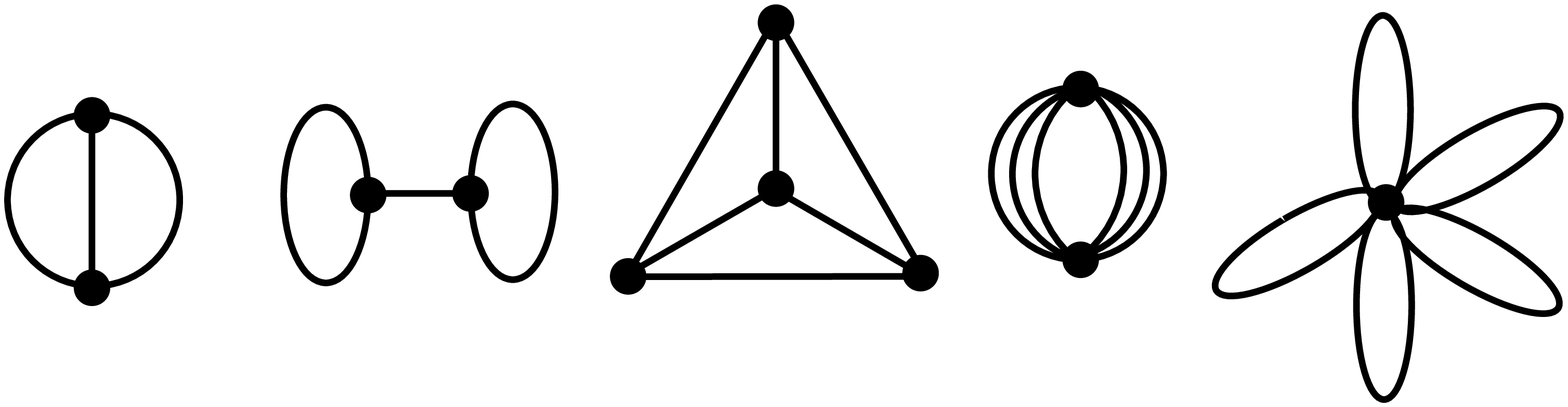}
\caption{From left to right, we have the abstract graph type of $\theta$-graphs, handcuff graphs, the tetrahedral graph, $\theta_n$-graphs, and bouquets.}
\label{fig:elem}
\end{figure}

\begin{figure}
\centering
\includegraphics[scale=0.3]{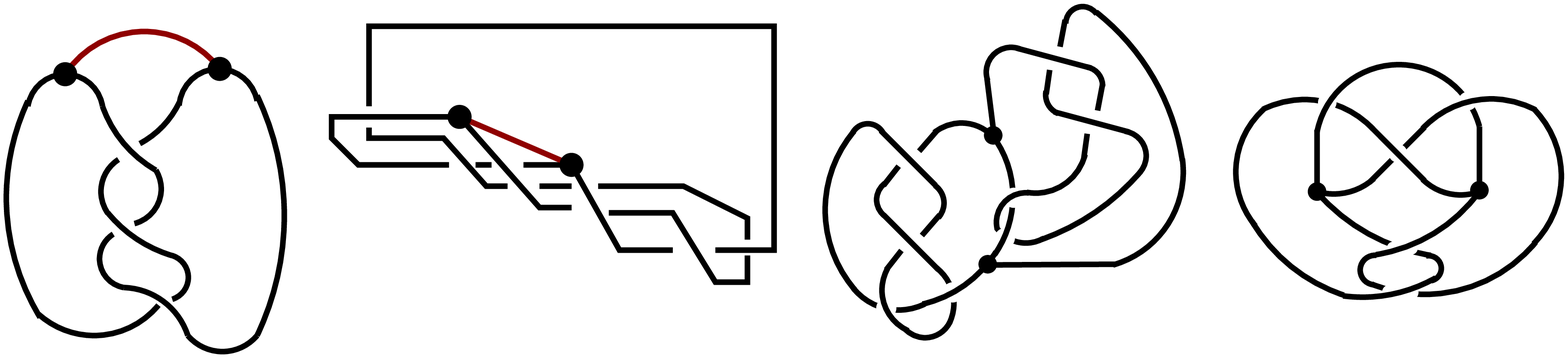}
\caption{Four examples of spatial $\theta$-graphs. From left to right: a trefoil knot with tunnel (in red); a knot with an arc (in red) marking the location of a crossing change; a $\theta$-graph whose every constituent knot is a figure 8 knot; the Kinoshita graph.}
\label{fig:thetagraphs}
\end{figure}

Some of the most important spatial graphs have very few vertices and edges, and are thus abstractly planar. Important classes include spatial $\theta$-graphs, handcuff graphs, the tetrahedral graph,  $\theta_n$-graphs, and bouquets. We depict the abstract graph type for these graphs in Figure \ref{fig:elem}. Of these families, spatial $\theta$ and $\theta_n$ graphs have received the most attention in the literature. One reason spatial $\theta$-graphs are so prevalent is that we can create one by attaching an arc to a knot so that the endpoints of the arc are distinct points on the knot.  This construction arises naturally in knot theory, where the arc may record some information about the knot $K$. Typical examples include knot tunnels or an arc recording the location of some crossing change, as in the first two diagrams of  Figure \ref{fig:thetagraphs}.  The knot $K$ becomes a cycle, or \defn{constituent knot}, in the resulting spatial $\theta$-graph. What can we say about the other constituent knots?  Perhaps suprisingly, Kinoshita \cite{Kinoshita87} showed that, given three knots, there is a spatial $\theta$-graph whose three constituent knots are precisely the three given knots.  An example of a  $\theta$-graph whose three constituent knots are all the Figure 8 knot is shown in Figure \ref{fig:thetagraphs}. Kinoshita's construction can be applied recursively to construct $\theta_n$-graphs whose constituent knots are specified beforehand. 

The trivial $\theta$-graph has every cycle an unknot. Are there other $\theta$-graphs with this property? A spatial graph having the property that no collection of disjoint cycles is a nontrivial link is a \defn{ravel}. Similarly, a spatial graph that has the property that every proper subgraph is trivial has the \defn{Brunnian property}. A knotted graph with the Brunnian property is \defn{Brunnian}, \defn{almost unknotted}, or \defn{minimally knotted}. A $\theta$-graph is Brunnian if and only if it is a ravel; but the same is not necessarily true for other graphs.  

Kinoshita \cite{Kinoshita1} also provided the first example of a Brunnian $\theta$-graph, now named after him. It is the rightmost diagram in Figure \ref{fig:thetagraphs}. Wolcott \cite{Wolcott} later generalized this construction to a family now known as the Kinoshita-Wolcott graphs. More examples of Brunnian $\theta$-graphs are given in \cite{Livingston} and \cite{JKLMTZ}. Suzuki \cite{Suzuki} generalized Kinoshita's construction to $\theta_n$ graphs.  Every abstractly planar graph without degree zero and degree one vertices has a Brunnian spatial embedding \cite{Kawauchi, Wu93}. Ravels are of interest to chemists \cite{Hyde}; Flapan and Miller \cite{FlapanMiller} have constructed many examples.

How can we be sure that Kinoshita's graph really is knotted, or indeed that any given spatial embedding of an abstractly planar graph really is knotted? 

An equivalence between spatial graphs takes the constituent knots of one to the constituent knots of the other (see \cite{Kauffman}). Thus, if one spatial graph has a constituent knot $K$ and another has no constituent knot of the same knot type, the graphs can't be equivalent. This doesn't help us show minimally knotted spatial graphs are knotted, though; we need other tools. As always in knot theory, we might ask for an invariant and there are some very nice invariants available. In general, Brunnian graphs and ravels provide good tests for the strength of invariants of spatial graphs. The three most popular are the Yamada polynomial \cite{Yamada}, Litherland's version of the Alexander polynomial \cite{Litherland}, and Thompson's polynomial invariant \cite{Thompson}. This last polynomial is defined recursively, but is zero if and only if the graph is trivial. It is based on an earlier algorithm of Scharlemann and Thompson \cite{ST-detecting} for determining if a spatial graph is unknotted. Their results were also adapted by Wu \cite{Wu92}, who showed that a spatial graph is unknotted if and only if each cycle bounds a disc disjoint from the rest of the graph.

We can also turn to other tools from topology and algebra. How accessible these are, depends, of course, on one's background and interests. Kinoshita and Suzuki used Alexander ideals to prove the nontriviality of their Brunnian graphs. However, McAtee, Silver, and Williams \cite{McAtee} point out that Suzuki's proof contains an error. The first complete proof of their nontriviality is likely given by Scharlemann \cite{Scharlemann}, using topological techniques stemming from the braid groups. In \cite[Example 22]{Ozawa}, the topology of surfaces containing the spatial graph is used to prove the Kinoshita graph is knotted, and in \cite{McAtee}, quandle colorings are used. Perhaps the simplest proofs that the Kinoshita graph is knotted rely on classical knot invariants. The article \cite{OT} provides two. In \cite{JKLMTZ}, a combination of handlebody theory and rational tangles are applied to an infinite family of $\theta$-graphs.   One popular and beautifully simple approach for $\theta_n$-graphs is to use branched covers. Livingston \cite{Livingston} uses these to prove Suzuki's graphs are nontrivial and Calcut and Metcalf-Burton \cite{Calcut} use them to show Kinoshita's graph is prime, in a sense which we now explore.

Whenever mathematicians are introduced to some new mathematical object, we want to be able to create more of them and to understand how the object fits into the larger context of known mathematics. For the remainder, we take up the question of creating new spatial graphs from old ones and understanding how spatial graphs are related to knot theory and 3-manifold theory.

\begin{figure}
\centering
\includegraphics[scale=0.3]{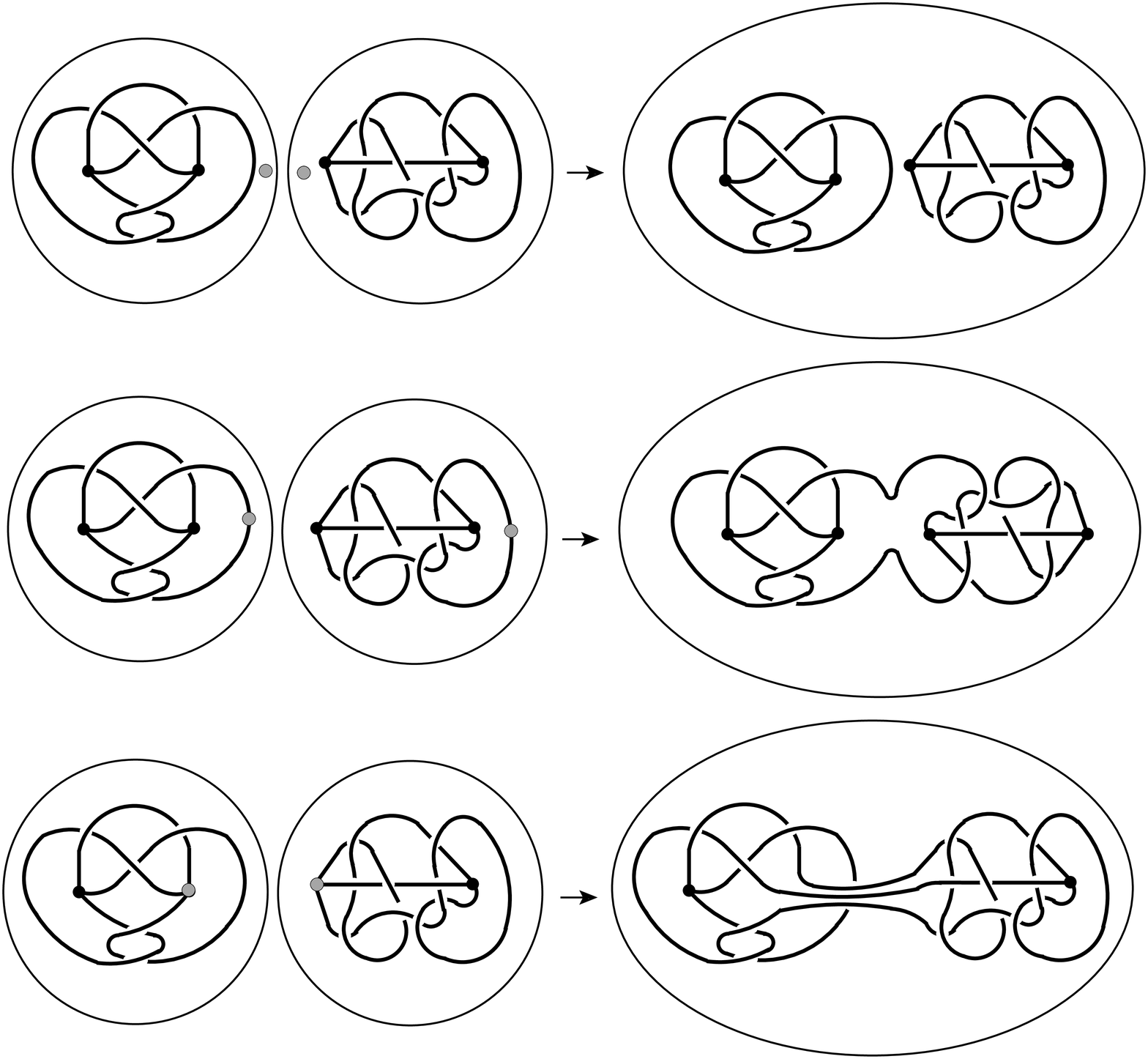}
\caption{From top to bottom we have the distant sum, the connected sum, and the trivalent vertex sum of two spatial $\theta$-graphs. In each case, the large circles and ellipses denote distinct copies of $S^3$ and the gray dots indicate the points where the summing occurs.}
\label{fig:summing}
\end{figure}

Throughout the study of manifolds, the connected sum is one of the most important methods of combining two manifolds. Recall that if $M_1$ and $M_2$ are manifolds of the same dimension, we form their connected sum $M_1 \# M_2$ by removing an open ball from each of them and gluing the resulting manifolds with boundary together along the new spherical boundary components. The centers of the balls are called the \defn{summing points}. Classical results show that if the summands are connected and we pay attention to orientations, the sum is unique up to homeomorphism. In particular $S^3 \# S^3$ is homeomorphic to $S^3$. When we consider 3-manifolds containing spatial graphs, on the other hand, there are  potentially many more options. For starters, we have a choice of where to perform the sum. Given two spatial graphs $G_1$ and $G_2$ in distinct copies $M_1$ and $M_2$ of $S^3$, we make the choice by picking summing points $p_1 \in M_1$ and $p_2 \in M_2$. We can pick the points to be both disjoint from the graphs, or we can pick them to be contained in the interiors of edges in the graphs, or we can pick them both to be vertices of the graphs having the same degree. We'll denote the result by $G_1 \#_k G_2$ where $k = 0$ if the points are disjoint from the graphs; $k=2$ if the points are interior to edges; and, otherwise, $k$ is the degree of the vertices\footnote{This creates conflicting notation when $k=2$, but we will ignore this.}. Figure \ref{fig:summing} depicts the the case when $k = 0$ (the \defn{distant sum}), $k=2$ (the \defn{connected sum}), and $k = 3$ (\defn{the trivalent vertex sum}). Even for a fixed $k$, in general, $G_1 \#_k G_2$ is not uniquely defined. Summing operations are associative.  If $G$ is equivalent to $G_1 \#_k G_2$ and neither is a trivial $\theta_k$-graph, then we say that $G$ is \defn{$k$-composite}. If $G$ is neither trivial nor $k$-composite, it is \defn{$k$-prime}.

For simplicity, let's consider only spatial $\theta$-graphs. We also assume our $\theta$-graphs are oriented. To orient a $\theta$-graph, choose one vertex as source, one vertex as sink, and color the edges red, blue, and green. We then restrict $\#_3$ so that, when forming $M_1 \#_3 M_2$,  the summing point $p_1$ is the sink vertex of $G_1$, the summing point $p_2$ is the source vertex of $G_2$, and the gluing map takes red, blue, and green endpoints to red, blue and green endpoints respectively.  Under the operation $\#_3$, the set $\mathbb{G}$ of oriented $\theta$-graphs in $S^3$ is particularly rich.  The operation $\#_3$ is well-defined \cite{Wolcott} and it makes the set $\mathbb{G}$ into a semigroup with the trivial graph as the identity. In other words, the operation $\#_3$ is associative and if $G \in \mathbb{G}$ and if $T$ is the trivial $\theta$-graph, then $G \#_3 T = T \#_3 G = G$. In general, elements of $\mathbb{G}$ do not have inverses. If they did, our semigroup would be a group. The center of the semigroup is the subset of elements commuting with all other elements. In our case, it consists of the $\theta$-graphs that are a connected sum of a trivial $\theta$-graph and a knot. Elements of $\mathbb{G}$ have 3-prime factorizations. This means that for each nontrivial $G \in \mathbb{G}$, there are 3-prime elements $G_1, \hdots, G_n$ such that $G = G_1 \#_3 \cdots \#_3 G_n$. Furthermore, except for the fact that elements of the center commute with all other elements,  this factorization is unique \cite{MT}. Conjecturally, (3-manifold, graph) pairs more generally have unique prime factorizations \cite{HAM}.

For $\theta$-graphs, the property of being Brunnian also persists under $\#_k$. Indeed, for $G_1, G_2 \in \mathbb{G}$, the sum, $G_1 \#_3 G_2$ is Brunnian if and only if $G_1$ and $G_2$ are both Brunnian. (Exercise!) For general spatial graphs, the property of being Brunnian may not persist under trivalent vertex sum. (Another exercise!) The property of being a ravel does persist under trivalent vertex sum. However, for $k \geq 4$, the property of being a ravel need not persist under $\#_k$. (Briefly: there exist knots with essential tangle decompositions and such knots result from summing bouquets.) Even for $\theta_k$ graphs, if we choose the gluing map for the connected sum to be very complicated, we may end up with knotted cycles after performing the sum.

\begin{figure}
\centering
\labellist
\small\hair 2pt
\pinlabel $e$ [l] at 72 183
\pinlabel $B$ [b] at 18 160
\pinlabel $B$ [br] at 453 196
\pinlabel $G$ [t] at 536 40
\pinlabel $\alpha$ [l] at 71 140
\endlabellist
\includegraphics[scale=0.3]{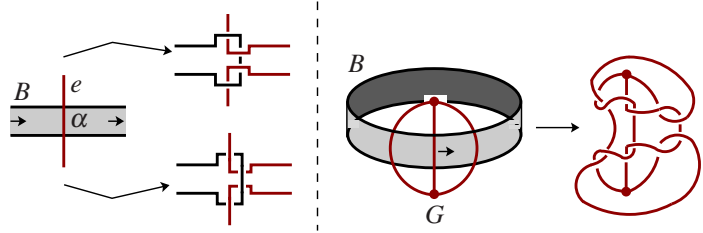}
\caption{The left side shows how to form $G(B)$. We take the union of the graph$G$ and the buckle $\boundary B$, but wherever an edge $e$ of $G$ intersects $B$ as on the left, we replace it with either of the pictured two ``belt buckles.'' On the right, we see that using the indicated belt $B$ to buckle the trivial theta-graph $G$ produces the trivalent vertex sum of the Kinoshita graph with its mirror image.}
\label{fig:buckle}
\end{figure}

How can we construct infinitely many \emph{3-prime} Brunnian $\theta$-graphs? The Kinoshita and Kinoshita-Wolcott graphs are 3-prime \cite{Litherland, Calcut}, as are the Brunnian $\theta$-graphs found in \cite{JKLMTZ} (see \cite{TT} for an indication of how this might be proved). Here is a very general method (essentially found in \cite{SimonWolcott}) that likely produces arbitrarily complicated Brunnian $\theta$-graphs, most of which are probably 3-prime.  For a spatial graph $G \subset S^3$, a new spatial graph $G(B)$, called a \defn{buckling} of $G$, is determined by a choice of oriented annulus $B = S^1 \times [0,1]$, called a \defn{belt}, intersecting $G$ in intervals of the form $\{\text{point}\} \times [0,1]$. We create $G(B)$ as follows. At each intersection arc $\alpha$ between $B$ and $G$ we replace a neighborhood of $\alpha$ in $S^3$ with a \defn{belt buckle} as on the left of  Figure \ref{fig:buckle} and include the remaining portions of $\boundary B$ as part of $G(B)$. The right of Figure \ref{fig:buckle} shows how a certain buckling of the trivial $\theta$-graph produces the trivalent vertex sum of the Kinoshita graph with its mirror image.

Not very much is known about how buckling affects a spatial graph. It is not difficult to see, however, that $G$ and $G(B)$ are abstractly isomorphic and that if $e$ is an edge of $G$ intersecting the band $B$, then the subgraphs $G(B) - e$ and $G - e$ are equivalent. In particular, if $G$ has the Brunnian property, and if $B$ intersects every edge of $G$, then $G(B)$ also has the Brunnian property. It seems difficult to determine whether or not $G(B)$ is trivial. Nevertheless, we conjecture:

\begin{conjecture}\label{conj:buckling nontrivial}
There is no $\theta$-graph $G$ in $S^3$ such that there is a belt $B$ intersecting all the edges of $G$ with $G(B)$ either the trivial $\theta$-graph or the Kinoshita graph.
\end{conjecture}

The exteriors of $G \in \mathbb{G}$ provide fertile ground for exploring topological and geometric questions. By work of Thurston (see Colin Adam's article in this encyclopedia), most knot complements admit a geometric structure based on 3-dimensional hyperbolic space. Likewise, the exteriors of most spatial graphs similarly admit hyperbolic structures. If our spatial graph is hyperbolic, we have access to a number of useful invariants. In particular, we may consider the \emph{volume} of the exterior. The simplest type of hyperbolicity for a spatial graph is \defn{hyperbolicity with parabolic meridians} \cite{Heard}. It follows from work of Thurston (see \cite[Corollary 2.5]{Heard}) that $G \in \mathbb{G}$ is hyperbolic with parabolic meridians if and only if it is $2$-prime and the exterior of $G$ does not contain an essential torus.  In particular, if $G \subset S^3$ is a Brunnian $\theta$-graph, then $G$ is hyperbolic with parabolic meridians whenever its exterior does not contain an essential torus. In \cite{JKLMTZ}, there is an example of a buckling producing a Brunnian $\theta$-graph having an essential torus in its exterior. That graph is, therefore, non-hyperbolic. Thurston's work can also be used to show that if $G = G_1 \#_3 G_2$ is some trivalent vertex sum of nontrivial elements of $\mathbb{G}$, then $G$ is hyperbolic with parabolic meridians if and only if both $G_1$ and $G_2$ are. This suggests that volume $\vol(G)$ is a particularly interesting invariant for $G \in \mathbb{G}$.  We ask (based on \cite{Adams}):

\begin{question}
Suppose that $G_i \in \mathbb{G}$ for $i = 1,2$ are both hyperbolic with parabolic meridians. How different can $\vol(G_1 \#_3 G_2)$ and $\vol(G_1) + \vol(G_2)$ be?
\end{question}

\begin{figure}
\centering
\includegraphics[scale=0.35]{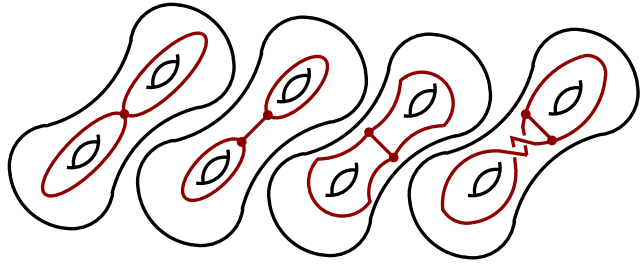}
\caption{Four different spines for a genus 2 handlebody}
\label{fig:handlebody}
\end{figure}

Finally, returning to topology, we consider the uniqueness up to homeomorphism of the exterior of connected spatial graphs. The relevant issues are illuminated by considering not just a spatial graph $G$ but also a regular neighborhood of it $N(G)$. The neighborhood $N(G)$ is a 3-manifold with boundary. It is called a \defn{handlebody} and the graph $G$ is called a \defn{spine} for the handlebody. The boundary of $N(G)$ is a connected, orientable surface. Its genus is called the \defn{genus} of the handlebody $N(G)$. Handlebodies may have many different spines. That is, if $N(G)$ is a handlebody with a spine $G$, there may exist many other graphs $G' \subset N(G)$ such that $N(G) = N(G')$.  Indeed, if the genus of the handlebody is at least 2, it will have infinitely many spines. Figure \ref{fig:handlebody} depicts four different spines for a genus 2 handlebody.  We say that two spatial graphs $G$ and $G'$ are \defn{neighborhood-equivalent} if $N(G)$ and $N(G')$ are ambient isotopic. Equivalently, $G$ and $G'$ are neighborhood-equivalent if there exists a handlebody in $S^3$ having two spines, one of which is equivalent to $G$ and the other to $G'$.  If two graphs are equivalent then they are neighborhood-equivalent, but the converse is not necessarily true. Since the exterior of a spatial graph $G$ is identical to the exterior of $N(G)$, neighborhood-equivalent spatial graphs have homeomorphic exteriors. In particular, spatial graphs are not determined by their complements, unlike knots in $S^3$ \cite{GordonLuecke}.

If two spatial graphs are neighborhood-equivalent, we can also ask how their constituent knots are related. For $\theta$-graphs, this question was studied extensively in \cite{Taylor1, Taylor2}, where it was connected to an operation on knots and 2-component links called \defn{boring}. Many operations in knot theory, such as rational tangle replacement on knots (an important operation in studying DNA, see e.g. \cite{ErnstSumners}),  are examples of boring. One result (see \cite[Theorem 6.5]{Taylor2}) arising from that work is that two Brunnian $\theta$-graphs are equivalent if and only if they are neighborhood-equivalent.  This suggests that Brunnian $\theta$-graphs may be determined by their complements.
\begin{conjecture}
If two Brunnian $\theta$-graphs have homeomorphic exteriors, then they are equivalent.
\end{conjecture}

In general, the topology of Brunnian graphs (not necessarily, $\theta$-graphs) is an area ripe for further study.

\section*{Acknowledgments} I would like to thank the referee and my students Tara Brownson and Qidong He for helpful comments. This research was partially supported by a research grant from Colby College.

\bibliography{Bib-short2}
\bibliographystyle{plain}

\end{document}